\documentclass{article}
\usepackage{amsmath}
\usepackage{amscd}
\usepackage{amsthm}
\usepackage{amssymb} 
\usepackage{mathtools}
\usepackage{tikz-cd}
\usepackage{latexsym}

\usepackage{euscript}
\usepackage{epsfig}

\usepackage{graphics}
\usepackage{array} 
\usepackage{enumerate}
 
\usepackage{pictexwd,dcpic}

\usepackage{tikz}
\usetikzlibrary{calc,intersections,through,backgrounds,arrows.meta}

\usepackage{color}

\usepackage{geometry}
\usepackage{hyperref}

\newcommand{\bel}[1]{\begin{equation}\label{#1}}

\newcommand{\be}{\begin{equation}}

\newcommand{\ba}{\begin{eqnarray}}
\newcommand{\ea}{\end{eqnarray}}
\newcommand{\rf}[1]{(\ref{#1})}
\newcommand{\bi}{\bibitem}
\newcommand{\qe}{\end{equation}}

\newcommand{\R}{{\mathbb R}}
\newcommand{\N}{{\mathbb N}}
\newcommand{\Z}{{\mathbb Z}}

\newcommand{\F}{\mathbb{F}}
\newcommand{\T}{\mathbb{T}}
\newcommand{\K}{\mathbb{K}}
\newcommand{\G}{\mathbb{G}}

\newcommand{\Hom}{\mathrm{Hom}}
\newcommand{\lth}{\mathrm{length}}

\newtheorem{thesis}{Thesis}
\newcommand{\btl}[1]{\begin{thesis}\label{#1}}
\newcommand{\et}{\end{thesis}}

\theoremstyle{theorem}

\newtheorem{satz}{Proposition}[section]

\theoremstyle{corollary}
\newtheorem{coro}{Corollary}[section]
\theoremstyle{lemma}
\newtheorem{lemma}{Lemma}[section]
\theoremstyle{definition}
\newtheorem{defi}{Definition}[section]
\theoremstyle{remark}

\theoremstyle{remark}
\newtheorem*{rem}{Remark}
\theoremstyle{remark}
\newtheorem{ex}{Example}

\title{Geometry and Dress groups with non-symmetric cost functions}
\author{Lukas Silvester Barth, Parvaneh Joharinad,\\  J\"urgen  Jost, Walter Wenzel}
\begin{document}
\maketitle

\begin{abstract}
    A metric relation by definition is symmetric. Since many data sets are non-symmetric, in this paper we develop a systematic theory of non-symmetric cost functions. Betweenness relations play an important role. We also introduce the notion of a Dress group in the non-symmetric setting and indicate a notion of curvature.
\end{abstract}
\section{Introduction}
Andreas Dress was a pioneer in developing profound new relations between algebra and geometry, and several of his constructions had an important impact for applications, for example the reconstruction of phylogenetic trees. In particular in \cite{D} (partly rediscovering earlier work of \cite{Isb}), he developed a general theory of tight spans (often called hyperconvex hulls in the literature and also in the present paper). From the mathematical side, this amounts to a systematic and penetrating study of metric spaces. A metric is a positive symmetric relation between the points of a set, and it satisfies the triangle inequality. Starting with the work of Hausdorff (see \cite{Haus}), it has become one of the central and most fertile notions of modern mathematics, and recently, it has also become immensely useful in machine learning (see for instance \cite{joharinadJost2023}). \\
In another line of research (though not unconnected to the previous work), in \cite{DW}, Andreas Dress and the fourth named
author have introduced the {\it Tutte group of a
matroid} that controls questions concerning representability; this is an abelian group
with generators and relations.\newline
Similarly, in \cite{JW}, we have introduced
for any set with a betweenness relation a
corresponding abelian group with generators
and relations, too, which we have named the
{\it Dress group}, in  view of the work
of Andreas Dress within the theory of
metric spaces just alluded to. Roughly speaking, the
Dress group plays a very similar role
for metric spaces as the Tutte group
does for matroids.\\

In this paper, hoping to preserve the spirit of the work and the thinking of Andreas Dress, we want to initiate a line of research that generalizes the aforementioned ideas to non-symmetric relations. That is, we want to abandon the symmetry requirement for a metric (and occasionally even its other properties), and see what kind of theories can emerge.\\

We shall thus systematically develop the geometric foundations for non-symmetric relations. These relations could express dissimilarities, transportation costs, transition probabilities, the difficulty  to rewrite one program into another with a universal Turing machine, or other possibly non-symmetric relations. There are natural links to partial orders, which we explore through the concept of betweenness relations. 
We shall also introduce and 
study a non commutative version of this
Dress group. Though this group is
formally still further away from the
-- commutative -- Tutte group of a
matroid, there are still very corresponding
facts, see Proposition \ref{W.8}:\newline
The non commutative Dress group contains
a distinguished subgroup ${\K}_0$ that is
analogous to the {\it inner Tutte group
of a matroid} as studied in \cite {DW}. We also briefly point out a relation with the concept of path homology of \cite{GLMY}.\\
We then turn to the concept of curvature. We find that  the general formulation of curvature in terms of ball intersection properties from \cite{joharinadJost2019} can be naturally extended to the non-symmetric case. Finally, we indicate what  the concept of the {\it tight span of a metric space} (hyperconvex hull) 
as studied by Andreas Dress in \cite{D} would look like in the non-symmetric case. 
For that purpose, we consider function pairs $(f,g)$.
More precisely, if $c$ as in \eqref{hyp1}
equals a metric $d$, that is, is symmetric, then those functions $f$,
for which the pairs $(f,f)$ satisfy
\eqref{hyp1}, define the tight span of the
metric space $(S,d)$, see also \cite{Isb}.
This tight span is, similarly to the convex
closure of a subset of a Euclidean Space, 
a topologically connected metric space,
which contains the given metric space as
an isometric substructure.\\

Naturally, the theory sketched here is by no means as complete and rich as those of metric spaces, that is, when the cost function is symmetric. Here, we have just defined some basic concepts that need to be explored in further research and applied to situations in machine learning where the relations between data points are non-symmetric.

\section{Non-symmetric cost functions and betweenness relations}\label{nonsymm}
\begin{defi}\label{costfunc}
 A \emph{cost function}  on a set $S$ is a relation $c:S\times S\to [0,\infty]$ 
 satisfying for all
$p,q,r\in S$
\ba
\label{1}
 c(p,q)&=&0 \text{ if and only if } p=q\\
\label{2}
c(p,r)&\le& c(p,q)+ c(q,r).
\ea
\end{defi}
In contrast to a metric, we do not require the symmetry
$c(p,q)=c(q,p)$. We may call $c(p,q)$ the \emph{cost} for getting from $p$ to
$q$. When $c(p,q)=\infty$, we may say that $q$ cannot be reached from $p$. But
even if $c(p,q)=\infty$, $c(q,p)$ may be finite, that is, we may be able to
reach $p$ from $q$. \rf{1} eliminates the need for the qualification \emph{pseudo}. There are
some relations for which even \rf{2} does not hold, like the Kullback-Leibler
divergence in information geometry, for which nevertheless interesting
properties can be derived, but here we want to keep \rf{2}. 

Some possible interpretations are that $c(p,q)$ measures the cost or the time
it takes to get from $p$ to $q$ in $S$. Our terminology below will draw upon
the latter interpretation. \\

\begin{defi}
\label{def:cost_morphisms_and_cost_space_category}
  Let $(S_1,c_1), (S_2,c_2)$ be sets with cost functions. A map $f:S_1\to S_2$
  is called a \emph{cost morphism} if it is \emph{cost nonincreasing}, that is, for
  all $p,q\in S_1$
  \bel{2aa}
  c_2(f(p),f(q))\le c_1(p,q).
  \qe
  $(S,c)$ is called a \textit{cost space}. The \emph{category of cost spaces} is the category in which the objects are cost spaces and the morphisms are cost morphisms.
\end{defi}

\begin{ex}\label{ex1}
An important example in the sequel is the unit interval $I=[0,1]$ with the cost
function
\bel{2a}
c(s,t)=\begin{cases} t-s \text{ if } s\le t\\
  \infty \text{ else}.
\end{cases}
\qe
\end{ex}
\begin{ex}\label{ex2}
Another important example obtains when we identify the endpoints $0$ and $1$
of $I$ to obtain a circle $S$ with the cost function
\ba
\label{2b}
c(0,1)&=&0\\
\text{else  }\qquad c(s,t)&=&\begin{cases} t-s &\text{ if } s\le t\\
  1+t-s &\text{ if }s>t .
\end{cases}
\ea
\end{ex}
\begin{ex}\label{ex3}
We shall also use  the example of the unit interval $I=[0,1]$ equipped with the
cost function
\begin{equation}
  \label{hyp2}
  c_1(p,q)=
  \begin{cases}
    q-p& \text{ if }p\le q\\
    2(p-q)& \text{ if }p\ge q.
  \end{cases}
\end{equation}
This construction will be taken up again in Section \ref{hyp}. See also Lemma \ref{monoid} in this context. 
\end{ex}
Looking at \rf{2}, a question is whether there exist $q$ different from $p,r$ with equality. In
order to explore this, it might be useful to recall the concept of a
\emph{betweenness} relation $b(p,q,r)$, even though we shall want to give up the symmetry
requirement of that notion.
\begin{defi}\label{def:between}
A \emph{betweenness relation} is  a three-point relation $b(p,q,r)$ that  satisfies
\ba
\label{3}
\text{ if } b(p,q,r)& & \text{ then } p,q,r \text{ are distinct}\\
\label{5}
\text{ if } b(p,q,r)& & \text{ then not } b(q,p,r)\\
\label{6}
\text{ if } b(p,q,r)\text{ and } b(p,r,s) & & \text{ then } b(p,q,s)\text { and } b(q,r,s)\\
\label{7}
\text{ if } b(p,q,s)\text { and } b(q,r,s) & & \text{ then } b(p,r,s)\text{ and } b(p,q,r).
\ea
\end{defi}
\begin{rem}
   \rf{6} and \rf{7} can perhaps be better remembered in the following shorthand notation
   \begin{eqnarray}
       \label{6a} 123\text{ and } 134 &\Longrightarrow &124\text{ and } 234\\
       \label{7a} 124\text{ and } 234 &\Longrightarrow  & 134\text{ and }123\ .
   \end{eqnarray}
\end{rem}
In Def.~\ref{def:between}, \rf{3} is for convenience. We point out that we do {\bf not} require  the symmetry condition that if $ b(p,q,r)$ then also $b(r,q,p)$ (in shorthand: $123 \Longrightarrow 321$)
that is usually required for a betweenness relation. A non-symmetric betweenness relation is for example natural 
  for analyzing (partially) ordered structures, with $b(p,q,r)$ if
 $p<q<r$. \\
 We also do {\bf not} require $123 \text{ and } 234 \Longrightarrow 124 \text{ or }134$. This is illustrated in the following diagram
 \begin{equation}\label{dig}
 \begin{tikzpicture}
     \fill (11,7) circle (3pt);
  \fill (11,1) circle (3pt);
   \fill (13,3) circle (3pt);
  \fill (13,5) circle (3pt);
  \coordinate [label=below:$p_1$] (x) at ($ (11,.8)$);
  \coordinate [label=right:$p_2$] (x) at ($ (13.2,3)$);
  \coordinate [label=right:$p_3$] (x) at ($ (13.2,5)$);
   \coordinate [label=above:$p_4$] (x) at ($ (11,7.2)$);
  \draw [-{Latex[length=4mm]}] (11,1) -- (11,7);
   \draw [-{Latex[length=4mm]}](11,1) -- (13,3);
   \draw [-{Latex[length=4mm]}](13,3) -- (13,5);
 \draw [-{Latex[length=4mm]}] (13,5) -- (11,7);
 \end{tikzpicture}
 \end{equation}
Here, one could for instance assume that the cost of each arrow is $1$; in the reverse direction, the cost could be some large number, say $10$. $p_2$ is between $p_1$ and $p_3$, and $p_3$ is between $p_2$ and $p_4$, but neither of them is between $p_1$ and $p_4$, as one can directly go from $p_1$ to $p_4$ with cost $1$. 
\medskip

Note also that every metric
$d: S \times S \rightarrow
{\R}^+ \cup \{0,\infty\}$ defines a
betweenness relation by writing $b(p,q,r)$
if and only if $p,q,r$ are distinct and satisfy

\medskip

\centerline{$d(p,q) + d(q,r) = d(p,r) < \infty$.}

\medskip

Clearly, \rf{3} and \rf{5} are fulfilled.\newline
Concerning the first conclusion in \rf{6}, suppose that $b(p,q,r)$ and
$b(p,r,s)$. Then we get

\medskip

\centerline{$d(p,s) \leq d(p,q) + d(q,s)
\leq d(p,q) + d(q,r) + d(r,s) = d(p,r) + d(r,s) = d(p,s)$}

\medskip

and, hence, $d(p,q) + d(q,s) = d(p,s) < \infty$ as claimed.

\medskip

To verify the second conclusion in \rf{6}, we obtain:

\medskip

\centerline{$d(q,s) \leq d(q,r) + d(r,s) =
d(p,r) - d(p,q) + d(p,s) - d(p,r) \leq d(q,s)$,}

\medskip

whence $d(q,s) = d(q,r) + d(r,s)$. -- Note that all of the distances that are involved here are finite.

\medskip

Similarly, one verifies \rf{7} -- via two
completely dual arguments.

\bigskip

\begin{ex}\label{paths}
    As already remarked, from a partial order $\le$, we can define a betweenness relation, with $b(p,q,r)$ if $p\neq q \neq r$ and $p\le q \le r$. This implies in particular to set systems, that is $\mathcal{S}\subset \mathcal{P}(V)$, where the latter is the power set of a set $V$. The subset relation in $\mathcal{S}$ then yields the partial order. In particular, if whenever $\rho \subset \sigma \subset V$ and $\sigma \in \mathcal{S}$, then also $\rho \in \mathcal{S}$, the set system defines a simplicial complex. Without that condition, we can consider $\mathcal{S}$ as a hypergraph. There is another way to view a simplicial complex $\mathcal{S}$ as a partially ordered set system, as pointed out in \cite{GLMY}. For a vertex set $V$, we consider a system of  finite paths $(v_{i_0},v_{i_1},\dots, v_{i_n})$. (There may be repetitions, that is, $v_{i_{j+1}}=v_{i_j}$ for some $j$, but for the theory developed in \cite{GLMY}, they turn out to be irrelevant.) When the system contains for any path also all its subpaths, we have a simplicial complex, and a  path is the ordered set of vertices of a simplex. Again, the subpath relation provides us with a partial order, hence a betweenness relation.
    \end{ex}
\bigskip

Here, we want to explore betweenness relations arising from cost functions $c$
satisfying \rf{1}, \rf{2}.
\begin{defi}\label{betdef3}
  Let $(S,c)$ be a set with a finite cost function. 
  $q\neq p,r$ then is between $p$ and $r$ if
  \bel{b1}
  c(p,r)=c(p,q)+c(q,r).
  \qe
\end{defi}

For the example $(I,c)$ of \rf{2a}, then $q$ is between $s$ and $t\in I$ if
\begin{equation*}
   s<q <t.
\end{equation*}
For the example of the circle $(S,c)$, see Example \ref{ex2}, $q$ is between $s$ and $t$ if the
three points are in cyclic order, 
\begin{equation*}
  s<q<t,\quad q<t<s\quad \text{ or } \quad t<s<q.
\end{equation*}
In fact, this looks more natural than the case where we take the usual distance function $d(.,.)$ on the unit circle and say that $q$ is between $p$ and $r$ if $d(p,q)>0, d(q,r)>0$ and $d(p,q)+d(q,r)=d(p,r)$. When $p$ and $r$ are antipodal, then any other point is between them for that latter relation, but if they are not, only the points on some segment of the circle are between them. This becomes even more drastic on spheres of dimension $\ge 2$, again equipped with the standard metric. When $p$ and $q$ are antipodal, again any other point is between them, but if they are not, only the interior points on the unique shortest geodesic arc are between them. 

\subsection{More general non-symmetric cost functions }\label{weakconditions}

Even without dropping assumption \ref{2}, one can weaken the definition of a cost function in two additional ways in order to connect it to other mathematical areas of interest.
\begin{itemize}
    \item Eliminating condition \rf{1}, allowing negative values of $c$.  
    \item Replacing the ``if and only if'' condition in \rf{1} by a simple ``if'' condition.
    \end{itemize}

With the first relaxation,  the triangle inequality still always ensures non-negativity of $c(p,p)$. However, non-negativity of $c(p,q)$ for $p\neq q$ becomes a non-trivial requirement because it does not follow automatically as in the case of a symmetric cost function, where we would have $2c(p,q)=c(p,q)+c(q,p)\ge c(p,p) \ge 0$. 

Cost functions for which the condition \rf{1} is weakened by the second option, generalize the well-established pseudo-metrics with some famous examples such as the space of metric spaces equipped with  Gromov-Hausdorff metric or the space of metric measure spaces equipped with Gromov-Wasserstein metric \cite{Sturm23, Memoli11}. Such cost functions are also related to metrics referred to as \textit{Lawvere metric spaces} in the literature. Lawvere observed in \cite{lawvere1973metric} that such spaces correspond to categories enriched in the monoidal poset $P=(([0,\infty],\ge),+)$, where the morphisms are $\ge$-relations and the tensor product is addition. To see this, note that the composition operation of morphisms in a $P$-enriched category $\mathbf{C}$ corresponds to the triangle inequality: The unique Hom-object of $p,q \in\mathbf{C}$ corresponds to the number $c(p,q)\in [0,\infty]$ and given $p,q,r\in\mathbf{C}$, composition $\circ_{p,q,r}:c(p,q)+c(q,r)\to c(p,r)$ is a morphism in the poset $([0,\infty],\ge)$, which exists iff $c(p,q)+c(q,r)\ge c(p,r)$. Furthermore, the axioms of a monoidally enriched category require composition to be unital, which implies $c(p,p)=0$. We provide a more detailed description in Section \ref{enrich}.\\
Starting from any bounded pairwise weight function  $w:G\times G\to [0,C]$, where $C \in\mathbb{R}$ is some constant, on a set $G$ with $w(g,g)=0$, one can always define a  Lawvere metric on $G$, as pointed out to us by Janis Keck. The proposed Lawvere metric is defined by $c(p,q):=w(p,q)+C(1-\delta_{pq})$, where $\delta$ is the \textit{Kronecker-Delta}. The definition implies that $c(p,p)=0$ and if $p,q,r$ are all different
\begin{equation}
c(p,r)+c(r,q) = 
    w(p,r)+w(r,q) + 2C \ge w(p,q)+C =c(p,q).
  \label{eq:voidTriangleInequality}
\end{equation}
However, in this Lawvere metric space, the triangle inequality is ``void'' in the sense that it does not enforce any constraints. It only holds because of a constant that is added to the weights. Every bounded pairwise weight can be made into a metric space in this way by adding a comparatively large constant  (e.g. maximum weight) to the weights to make the cost of going from $p$ to $r$ and then proceeding to $q$ big enough that it always exceeds the cost of going directly from $p$ to $q$. A natural question is how one might filter out those metric spaces that are coming from pairwise weights that do not fulfill the triangle inequality.\\
Below we provide one possible resolution by defining categories of cost spaces in which the axioms of a cost function or Lawvere metric space hold only up to a constant. To make this precise, let us introduce the $\mathcal{O}$-notation. Let $f,g$ and $h$ be functions with domain $D$ and values in the  real numbers (possibly extended by $\infty$). Then we define
\begin{equation}
  \label{eq:O-notation}
  \begin{split}
    f \le g + \mathcal{O}(h)\quad\text{iff}\quad 
    \exists C<\infty ~:~f(x)\le g(x)+Ch(x)\quad\forall x \in D.
  \end{split}
\end{equation}
We also introduce a useful variation for bounded functions with two arguments, which is supposed to handle the problem described around equation \eqref{eq:voidTriangleInequality}.  
Let $f$ be any function from $D\times D$ to $\mathbb{R}$ and define $B(f):=\inf_{y\ne z}f(y,z)$ as well as $f^B(x,x'):=f(x,x')-B(f)$. Similarly, let $\{g_i\}_{i\in I}$ be functions with domain $D\times D$. Then we define $\mathcal{B}(1)$  as follows:
\begin{equation}
  \begin{split}
    f(x,x') \le \sum_{i\in I}g_i(x_i,x_i') + \mathcal{B}(1)\text{ iff }
    f^B(x,x') \le \sum_{i\in I}g_i^B(x_i,x_i')
  \end{split}
  \label{eq:B1inequality}
\end{equation}
for all $x,x',x_i,x_i'$.
We can use the above to define more general costs.
\begin{defi}\label{asymptoticCosts}
  An \emph{asymptotic cost function} or $\mathcal{O}$-\textit{cost function} on a set $S$ is a relation $c:S\times S\to [0,\infty]$ satisfying for all
 $p,q,r\in S$
 \ba
 \label{o1}
 c(p,p)&=&0 + \mathcal{O}(1)\\
 \label{o2}
 c(p,q)&\le& c(p,r)+ c(r,q) + \mathcal{O}(1).
 \ea
 \end{defi}
 Similarly, a $\mathcal{B}$-\textit{cost} is obtained by replacing $\mathcal{O}(1)$ 
 in \rf{o1} and \rf{o2} by $\mathcal{B}(1)$. The category of such $\mathcal{B}$-cost spaces then no longer contains metric spaces with a void triangle inequality as in eq.~\eqref{eq:voidTriangleInequality} and is more natural in this sense. 

Besides the exclusion of pathological spaces, this definition allows us to connect the theory of cost functions to the theory of formal languages, rewrite theory, Kolmogorov complexity and related ideas in computational linguistics, which we shall explain below.

\subsubsection{Conditional Kolmogorov complexity cost function}

Kolmogorov introduced the complexity named after him in \cite{kolmogorov1963complexity} and modern references include \cite{livitanyi2008KolmogorovComplexity} and \cite{shenUspensky2017}. We briefly introduce the elementary definitions.
Let $A$ be a set, let $A^*$ be the set of all finite sequences of elements of $A$, let $p,q\in A^*$, let $pq$ denote the concatenation of $p$ and $q$, let $\ell(p)$ denote the length (usually in bits) of $p$ and let $U$ be a universal Turing machine. We also restrict the codes on which $U$ operates to be prefix codes (meaning that no code is the prefix of another one) because this simplifies some expressions in the sequel. The \textit{conditional (prefix) Kolmogorov complexity} $K_U(p|q)$ is then defined as follows:
\begin{equation}
  \begin{split}
      K_U(p|q) := \min_{r}\{~\ell(r)~|~U(rq)=p~\}.
  \end{split}
\end{equation}
The intuition is: $K_U(p|q)$ is the length of (one of the) shortest program(s) $r$, running on $U$, that generates $p$ when given $q$ as input. A special case is $K_U(p):=K_U(p|\epsilon)$, where $\epsilon \in A^*$ is the empty sequence. $K_U(p)$ is also simply called the Kolmogorov complexity of $p$ and denotes the length of the shortest program that can generate the sequence $p$ on $U$.

$K_U$ depends on the universal machine $U$. However, by definition of a universal machine, $U$ can emulate any other machine. This means that, for any other universal machine $U'$, there is a (finite) sequence $r\in A^*$ such that, for all $p\in A^*$, we have $U(rp)=U'(p)$. One can think of $r$ as the compiler of $U'$ in the language (presented by) $U$. Since $\ell(r)<\infty$, this implies that, up to a constant that does not depend on $p$, $K_U(p)$ is well-defined, independently of $U$, and we therefore sometimes omit the index $U$. This is one of the reasons why one might want to speak about the value of functions up to some constant as in eq.~\eqref{eq:O-notation}. The relationship to geometric notions is established by the following proposition, taken from \cite[Theorem II.1]{KolmogorovTriangleInequality}:
\begin{satz}\label{triangleInequalityK}
  $c(p,q):=K_U(p|q, K_U(q))$ fulfills \rf{o2}, i.e.
  \begin{equation}
    \begin{split}
      K_U(p|q,K_U(q)) \le K_U(p|r,K_U(r)) + K_U(r|q,K_U(q)) + \mathcal{O}(1).
    \end{split}
  \end{equation}
  \label{eq:KolmTriangleInequality}
\end{satz}
Furthermore, note that $K_U(p|q)$ is not symmetric because if $q$ is a sequence that contains $p$ as a prequel, then $K_U(p|q)$ is small because one only has to forget part of $q$ but $K_U(q|p)$ might be large if $q$ is much longer than $p$. Finally, also $K_U(p,p)$ is usually only $0$ up to a constant because to generate $p$ usually requires some non-empty program that contains the code ``return $p$'' or  similar. Hence, to make $(p,q)\mapsto K_U(p|q,K_U(q))$ into some kind of cost function, we really need to weaken the axioms to those specified in Definition \rf{asymptoticCosts}. 

Since $K_U(p|q)$ can also be understood as the shortest possible rewrite of $q$ into $p$, eq.~\eqref{eq:KolmTriangleInequality} also provides a link to rewrite theory and thus formal languages and their grammar.

\section{Swiftest curves}\label{swift}
We shall now extend the notion of a (shortest) geodesic developed in \cite{JW}
to non-symmetric relations. We shall give the concepts new names, again
derived from classical Greek.\\

We consider cost functions as described in Definition \ref{costfunc}.
\begin{defi}\label{tach}
 For distinct $p,q \in S$, a  \emph{tachistic} (from Greek
 \emph{tachistos} = swiftest) from $p$ to $q$ is a map
$g: J \rightarrow S$ defined on a subset $J$ of an interval $[a,b]$ in ${\R}$
with $a,b \in J$ that satisfies  $g(a) = p$ and $g(b) = q$ and for all $s<t,
s,t \in J$
\bel{33}
c(g(s),g(t))=t-s
\qe
and that cannot be  extended 
to some larger subset of $[a,b]$ as a map with values in $S$ satisfying this
property. 
\end{defi}
It follows from Zorn's Lemma that such   tachistics always exist.
\begin{defi}\label{general}
 For $p,q \in S$, a \emph{chronodesic} (from Greek \emph{chronos} =
  time) from $p$ to $q$ is a
map $g: J \rightarrow S$ defined on a subset $J$ of a compact  interval
$[a,b] \subseteq \R$ with $a,b \in J$ and satisfying the following conditions:
 \begin{itemize}
\item[(C1)] $g(a) = p, \, g(b) = q$, and either
  \begin{itemize}
  \item[(C1a)]  $J=\{a,b\}$ and $c(g(a),g(b))=b-a$\\
    or
  \item[(C1b)]  we can find $t_0=a<t_1<\dots < t_n=b
    \in J$ for some $n\ge 2$ with $c(g(t_{i-1}),g(t_{i}))+c(g(t_{i}),g(t_{i+1}))
    =c(g(t_{i-1}),g(t_{i+1}))\\ =|t_{i+1}-t_{i-1}|$
  for $i=1,\dots , n-1$\\
(that means, $g$ is a tachistic map on any interval $[t_{i-1},t_{i+1}]$, and
$g(t_i)$ is between $g(t_{i-1})$ and $g(t_{i+1})$).
 \end{itemize}
\medskip
\item[(C2)] There does not exist a continuation
$\tilde{g}: \tilde{J} \rightarrow S$ of $g$ to some set $\tilde{J}$ with
$J \subsetneq \tilde{J} \subseteq [a,b]$ with the property (C1b)  into the
same space $S$.
 \end{itemize}
   \end{defi}

In the same way that geodesic curves in Riemannian or metric geometric need not be shortest connections between their endpoints, also in non-commutative situations, chronodesics need not be tachistic. For instance, in diagram \eqref{dig}, the path from $p_1$ to $p_4$ via $p_2$ and $p_3$ is not tachistic, because one can get from $p_1$ to $p_4$ with lower cost directly.

\section{Dress groups}
  
From a betweenness relation, one can define the \emph{Dress group} \cite{JW}.
\begin{defi}
   For a betweenness relation $b$ on the nonempty set $S$,
the group $\tilde{\F} = \tilde{\F}_S$ is the free group generated by all
symbols $X_{p,q}$ for $p,q \in S$ with $p \neq q$.
Let $\tilde{\K} = \tilde{\K}_S$ denote the smallest {\it normal} subgroup
of $\tilde{\F}_S$ that contains all elements

\medskip

\centerline{$X_{p,r} \cdot X_{q,r}^{-1} \cdot X_{p,q}^{-1}$ whenever $b(p,q,r)$.}

\medskip

Then the -- in general -- {\it noncommutative Dress group}
$\tilde{\T} = \tilde{\T}_S$ is defined as the factor group

\medskip

\centerline{$\tilde{\T}_S = \tilde{\F}_S / \tilde{\K}_S$.}
\end{defi}
We also denote the image of $X_{p,q}$ in this quotient by
$\tilde{T}_{p,q}$. 

A direct consequence is 
\begin{coro}
  If $c: S \times S \rightarrow {\R^+} \cup \{0\}$
is a cost function, then we have a well defined homomorphism
$f: \tilde{\T}_S \rightarrow \R$ given by

\medskip

\centerline{$f(\tilde{T}_{p,q}) := c(p,q)$, where $b(p,q,r)$ holds if
and only if $c(p,r) = c(p,q) + c(q,r)$.}

\end{coro}\qed
\begin{rem}
    When we also want to consider cost functions that can become infinite, we should restrict homomorphisms to the subgroup generated by the images $\tilde{T}_{p,q}$ of those $X_{p,q}$ for which the cost $c(p,q)$ is finite. 
\end{rem}

Next we look at some generalization of {\it metric embeddings}.

\medskip

\begin{defi}Suppose that $S_1$ and $S_2$ are sets with
betweenness relations $b_1$ and $b_2$, respectively.
Then a {\it morphism} $\varphi: (S_1,b_1) \rightarrow (S_2,b_2)$
is a map $\varphi: S_1 \rightarrow S_2$ satisfying

\medskip

\centerline{whenever $b_1(p,q,r),$ then
$b_2(\varphi(p),\varphi(q),\varphi(r))$.}
\end{defi}

From the definitions, we conclude at once
\begin{satz}\label{W.2}
  Whenever $\varphi: (S_1,b_1) \rightarrow (S_2,b_2)$
is a morphism between sets with bet\-weenness relations, then $\varphi$
induces a canonical homomorphism
$\psi: \tilde{\T}_{S_1} \rightarrow \tilde{\T}_{S_2}$.
\end{satz}
\qed\\

For any nonempty set $S$ with a betweenness relation $b$, we have of
course a canonical epimorphism $\pi$ from the group $\tilde{\T}_S$
onto the commutative group ${\T}_S$ as studied in
\cite{JW}.
\newline
If, in particular, $S$ is a metric space with a metric $d$, then by
Proposition 3.6 in \cite{JW} (which is strongly related to Proposition \ref{W.2}
above) we have also a well defined homomorphism $f$
from ${\T}_S$ into the field of real numbers given by
$f(T_{a,b}) := d(a,b)$.
\newline
By combining these last two mentioned facts, we obtain at once

\medskip

\begin{satz} For any metric space $(S,d)$ we get a composed
homomorphism $f \circ \pi$:

\medskip

\centerline{$\tilde{\T}_S \rightarrow {\T}_S \rightarrow {\R}.$}

\end{satz}

\begin{rem} Note that ${\T}_S$ is -- in general -- {\it not}
the abelianization one gets by starting from $\tilde{\T}_S$.\newline
If, for instance, $S = \{s,t\}$ contains only two -- different --
elements, then there does not exist any betweenness relation
between these elements, whence $\tilde{\T}_S$ is -- by definition --
the free group generated by $2$ elements. --Its abelianization
is isomorphic to ${\Z}^2$. However, since $T_{s,t} = T_{t,s}$ holds
in the commutative Dress group ${\T}_S$, this group is the infinite
cyclic group -- and, hence, isomorphic to ${\Z}$.
\end{rem}
\medskip

\noindent
{\bf Conventions:} Suppose that $S$ is a -- finite or infinite --
set with at least $2$ elements and a betweenness relation $b$.\newline
Then ${\Z}^S$ is -- as usual -- the set of {\it all} maps from $S$
into ${\Z}$, while we denote by ${\Z}^S_{fin}$ the set all of these maps $f$
of {\it finite support}; that means, there are only finitely many elements
$s \in S$ with $f(s) \neq 0$.\newline
Moreover, for $s \in S$, the map ${\delta}_s$ signifies the map
with ${\delta}_s(s) = 1$ and ${\delta}_s(t) = 0$ for all $t \neq s$.
Then, all of these maps ${\delta}_s$ build a base of the free
${\Z}-$ module ${\Z}^S_{fin}$.\newline
Moreover, put
\begin{equation}
    {\G}_0 := \{\sum_{s \in S} \, n_s \cdot {\delta}_s
\in {\Z}^S_{fin} \, | \, \sum_{s \in S} \, n_s = 0\}.
\end{equation}

We can now prove the following

\medskip
\begin{satz}\label{W.8}
    Assume that the set $S$ has at least $2$
elements and that $S$ is equipped with a betweenness relation $b$.
Then we have a well defined homomorphism
$\psi: \tilde{\T}_S \rightarrow {\G}_0$ given by
\begin{equation}
    \psi(\tilde{T}_{s,t}) := {\delta}_t - {\delta}_s.
\end{equation}

Moreover, $\psi$ is surjective -- and, hence, an epimorphism.
If we thus denote by ${\K}_0$ the kernel of $\psi$, we get
\begin{equation}
    \tilde{\T}_S / {\K}_0 \, \simeq \, {\G}_0.
\end{equation}
\end{satz}
\begin{proof}
Clearly, the images of the homomorphism $\psi$ lie in ${\G}_0$,
whence $\psi$ is well defined -- by the definition of $\tilde{\T}_S$.
The only remaining nontrivial fact is that $\psi$ is surjective.
\newline
Suppose that $g := \sum_{s \in S} \, n_s \cdot {\delta}_s$ lies in ${\G}_0$.
We must prove that $g \in \psi(\tilde{\T}_S)$.\newline
We proceed by induction on $N := \sum_{s \in S} \, |n_s|$ (or, if one formally
prefers, by the half of this sum $N$).\newline
If $N = 0$, then $g$ is the neutral element in ${\G}_0$, whence
$g = \psi(1)$.\newline
Now assume that $N > 0$. Then, in view of $g \in {\G}_0$, there exist
elements $s,t \in S$ with $n_s < 0$ and $n_t > 0$.\newline
Put $g' := g + {\delta}_s - {\delta}_t$. By the induction hypothesis,
there exists an element $\tilde{T}' \in \tilde{\T}_S$ with
$\psi(\tilde{T}') = g'$. If we now put
$\tilde{T} := \tilde{T}' \cdot \tilde{T}_{s,t}$, we get
$\psi(\tilde{T}) = g$ as claimed.
\end{proof}
\bigskip
\begin{rem}
     Clearly, the ${\Z}-$module ${\G}_0$ does not
depend on the betweenness relation $b$, while the kernel ${\K}_0$
heavily depends on $b$. For any three pairwise distinct elements
$p,q,r \in S$ one has
\begin{equation}
    \tilde{T} :=
\tilde{T}_{p,r} \cdot \tilde{T}_{q,r}^{-1} \cdot \tilde{T}_{p,q}^{-1} \in {\K}_0\ ,
\end{equation}
but this product is -- in general -- only the neutral element if $b(p,q,r)$
holds.\newline
Proposition \ref{W.8} suggests to study ${\K}_0$ exhaustively.

\medskip

In \cite{DW}, we have studied a conspicuous subgroup of the Tutte group of a matroid -- called
the {\it inner Tutte group} -- which has properties
very similar to those of ${\K}_0$.
\end{rem}

\bigskip

\begin{lemma}\label{embedding}
Suppose that $X_1,...,X_n$ are certain -- pairwise distinct --
indeterminates and that $(X_i)_{i \in I}$ are further pairwise
distinct indeterminates, different from $X_1,...,X_n$.\newline
Consider the free group $G$ generated by all $X_1,...,X_n$
as well as the group $H$ generated by all $X_k, 1 \leq k \leq n$
and all $X_i, i \in I,$ and certain relations $X_i = f_i(X_1,...,X_n)$,
where each $f_i(X_1,...,X_n)$ is a product of certain powers
of the elements $X_1,...,X_n$ -- with positive or negative
exponents. Then the groups $G$ and $H$ are isomorphic.
More precisely, two inverse isomorphisms $g: G \rightarrow H$
and $h: H \rightarrow G$ are given by

\medskip

\centerline{$g(X_k) := X_k$ for $1 \leq k \leq n$,}

\medskip

\centerline{$h(X_k) := X_k$ for $1 \leq k \leq n, \,
h(X_i) := f_i(X_1,...,X_n)$ for $i \in I$.}

\end{lemma}
\begin{proof} Clearly, by the definitions, we have
$h(g(X_k)) = X_k$ and $g(h(X_k)) = X_k$ for all $k$ with $1 \leq k \leq n$.
\newline
Moreover, for $i \in I$ we have 

\medskip

\centerline{$g(h(X_i)) = (g \circ f_i)(X_1,...,X_n) = f_i(X_1,...,X_n) = X_i$}

\medskip

as claimed, where the second equation holds by the assumptions about the
functions $f_i$ -- and the definition of $g$.
\end{proof}

\medskip

\begin{ex} \label{ex4}
This example is closely related to Example \ref{ex2}; but, now, we consider "only" a discrete subset $S_n$ of the circle $S$:\newline
Suppose that the natural number $n$ satisfies
$n \geq 2$, and consider the roots of unity

\medskip

\centerline{$w_k :=
\exp(2 \pi i \cdot \frac{k}{n})$ for
$1 \leq k \leq n$,}

\medskip

as well as $S_n := \{w_1,...,w_n\}$.\newline
Moreover, $b(w_k,w_l,w_m)$ holds for pairwise different $k,l,m$ if and only if the counterclockwise arc
from $w_k$ to $w_m$ runs through $w_l$.
\newline
We can now apply Lemma \ref{embedding}
to the indeterminates
$X_k = \tilde{T}_{w_k,w_{k+1}}, k$ mod $n$, and the relations

\medskip

\centerline{$\tilde{T}_{w_k,w_{k+d}} =
{\Pi}_{1 \leq j \leq d} \,\, 
\tilde{T}_{w_{k+j-1},w_{k+j}}$}

\medskip

for $2 \leq d \leq n-1$ and $k+d$ mod $n$.\newline
We conclude that $\tilde{\T}_{S_n}$ is (isomorphic to) the free group generated by $n$ elements.

\end{ex}

\bigskip

\begin{ex} \label{directedgraph}
Suppose that $G = (V,E)$ is a finite {\it directed} graph; this means, that $E$ consists of pairs $(v,w)$ for distinct vertices $v,w \in V$. Assume furthermore that for any two distinct vertices
$u,v \in E$ there exists at most one path in $G$
with starting point $u$ and endpoint $v$; this
means in particular that $G$ does not contain any
directed circular path.\newline
Consider the -- natural -- betweenness relation
$b$ on $V$ by writing $b(u,v,w)$ for any three
pairwise distinct vertices $u,v,w$ if and only if
$v$ lies on some -- in this case unique -- path
from $u$ to $w$.\newline
Then $\tilde{T}_V$ is (isomorphic to) the free
group generated by all elements
$\tilde{T}_{v,w}$ for which either $(v,w)$ is
an edge in $E$ or $w$ is unreachable from $v$.

\medskip

Namely,  exactly those pairs $(v,w)$ which
are not listed here have the property that there
exists a -- by assumption unique -- path $(v_0,...,v_l)$ from $v$ to $w$ of length
$l \geq 2$. This means that

\medskip

\centerline{$\tilde{T}_{v,w} =
{\Pi}_{1 \leq j \leq l} \,\, 
\tilde{T}_{v_{j-1},v_j}$.}

\medskip

Hence, Lemma \ref{embedding} yields what we want.

\end{ex}

\begin{ex}\label{pathhomology}
 We can use the preceding example to compare the construction of the Dress group with that of the path homology groups of \cite{GLMY}. There, one considers a system $\mathcal{S}$ of finite paths  $(v_{i_0},v_{i_1},\dots, v_{i_n})$ with the $v_{i_j}\in V$, and assumes that every subpath of a path in $\mathcal{S}$ is also in $\mathcal{S}$. One can then define a boundary operator
 \begin{equation}
    \partial  (v_{i_0},v_{i_1},\dots, v_{i_n}) =\sum_j (-1)^j (v_{i_0},\dots, \widehat{v_{i_j}},\dots ,v_{i_n})
 \end{equation} which squares to 0,
 \begin{equation}
     \partial \circ \partial =0,
 \end{equation}
 and path homology groups. As explained in Example \ref{paths}, this generalizes the homology theory of simplicial complexes. So, here, one gets even a family of groups, as many as the cardinality of $V$. These are thus different from the single Dress group that we construct, but the Dress group can be constructed in a much wider setting than these path homology groups. In particular, as follows again from that example, the Dress group is defined for hypergraphs (without having to embed them into simplicial complexes, as in \cite{Ren}), and not only for simplicial complexes.
\end{ex}

\bigskip

In any case, from our perspective, we may ask whether on some metric space
$(S,d)$, for three or more distinct points $p_1,\dots ,p_n$, there is some $q$ that is
between any pair $(p_i,p_j)$ for $i\neq j$, or, more precisely, to quantify
the deviation, that is, to which extent the best choice $q$ violates equality
in \rf{2} for all such pairs.

A useful aspect is that with a group structure, one can derive further
relations by algebraic computations.

\section{Pretopologies}\label{pretop}
We recall the concept of a \emph{pretopological space}, see \cite{J}, also called  a \textit{\v{C}ech closure space} in the literature, after \cite{Cech}.
\begin{defi}
A set $X$ with power set $\mathcal{P}(X)$ is  a \emph{pretopological space} if it possesses a preclosure
 operator   $~ \overline{\bullet}\ $ with the following properties
 \begin{enumerate}
 \item[(i)] $\overline{\emptyset}=\emptyset$.
 \item[(ii)] $A\subset \overline{A}$ for all $A\in \mathcal{P}(X)$.
 \item[(iii)] $\overline{A \cup B}=\overline{A}\cup \overline{B}$ for
   all $A,B \in \mathcal{P}(X)$.
 \end{enumerate}
 $A\in \mathcal{P}(X)$ is called \emph{closed} if $A= \overline{A}$.
\end{defi}
We also recall that such an $X$  is  a topological space iff the preclosure operator in addition
 satisfies
 \begin{enumerate}
 \item[(iv)] $\overline{\overline{A}}=\overline{A}$ for all $A \in \mathcal{P}(X)$.
 \end{enumerate}

 In our context, $\overline{A}$ may be interpreted as  that part of $X$ that you
can reach from $A$ by applying some operation. By (i) then nothing can be
reached from nothing. By (ii), all starting points can be reached. By (iii),  from a union of
starting sets nothing  more can be reached than the combination of what 
can be reached from each single set.  In a  directed graph $\Gamma$,
one can define the preclosure of a set of vertices as the union of this
set with the set of all forward 
 neighbors of these vertices.  Conversely,
from a pretopological space, we can
construct a directed graph by connecting each $x$ with all the other 
elements of $\overline{\{x\}}$.\\
Another  example of a preclosure operator arises from a
dynamical system. For concreteness
\ba
\label{top-ds1}
\dot{x}(t)&=&F(x(t)) \text{ for } x\in \R^d,t>0\\
\label{top-ds2}                 
x(0)&=&x_0
\ea
for some uniformly Lipschitz continuous  $F$, which then possesses a 
unique 
solution  with initial values \rf{top-ds2} for all
$t\ge 0$.  For $A\subset \R^d$, we  put
\bel{top-ds3}
\overline{A}^T:=\{ x(t), 0\le t \le T\} \text{ where } x(t) \text{ is
  a solution of \rf{top-ds1} with } x(0) \in A.
\qe
For each $T>0$, this then defines a preclosure operator. The closed sets
 are the forward invariant sets of  \rf{top-ds1}.\\

We can also consider a collection $X$ of programs, and let
$\overline{A}$ contain the programs that can be reached from those in $A$ by a
predefined number of steps. Depending on which types of concatenation of
programs we allow, however, (iii) need not be satisfied. An example is
genetic recombination of strings. From a single string, one can reach nothing
else by recombination, but from recombining two different ones, one may reach
many others. Thus, we may have to work in a context somewhat more general than
that of pretopological spaces. For a systematic treatment of various notions
of such structures generalizing or extending that of a pretopology, see \cite{Sta1}.\\

\begin{rem}
It is also possible to lift the construction to a categorical setting. We provide a definition that is a slight generalization of the usual notion of a ``universal closure operation'' as given, for example, in \cite[Section 5.7]{Bor}.
\end{rem}
\begin{defi}\label{preclos}
Let $\mathcal{C}$ be a finitely complete category. A \textit{universal preclosure operation} on $\mathcal{C}$ consists in assigning, for every subobject $S \xhookrightarrow{} C$ (also written $S\subset C$) in $\mathcal{C}$, another subobject $\overline{S} \xhookrightarrow{} C$ called the \textit{preclosure of $S$ in $C$}, subject to the following conditions:  
 \begin{enumerate}
 \item[(i)] $S\subset \overline{S}$ for all subobjects $S$ of $C$.
 \item[(ii)] $S\subset T \Rightarrow \overline{S}\subset \overline{T}$ for
   all subobjects $S,T$ of $C$.
   \item[(iii)] $f^{-1}(\overline{S})= \overline{f^{-1}(S)}$ whenever $f:B\to C$ is a morphism in $\mathcal{C}$.
 \end{enumerate}
\end{defi}
 In the case of a pretopological space, the axioms require $\overline{A \cup B}=\overline{A}\cup \overline{B}$ for $A,B \in \mathcal{P}(X)$ instead of (ii). However, since this implies $A\subset B \Rightarrow \overline{A}\subset \overline{B}$, the above definition generalizes the definition of the pretopological space, while we recover a usual universal closure operation in the sense of \cite{Bor} upon the imposition of idempotency.

\bigskip

\begin{defi}
\label{def:continuity}
A map $f:Z\to X$ between pretopological spaces is \emph{continuous} iff
\bel{top11}
f(\overline{B})\subset \overline{f(B)}
\qe
for any subset $B$ of $Z$.
\end{defi}
More generally, this is also meaningful when $X$ possesses a preclosure operator
that only satisfies (i) and (ii), but not necessarily (iii) in Definition \ref{preclos}. When the pretopology satisfies condition (iv), Definition \ref{def:continuity} properly specializes to the usual notion of continuous map between topological spaces \cite[Lemma 4.1.5 and 4.1.6]{J}. 
\begin{defi}
The \emph{category of pretopological spaces} is the category whose objects are pretopological spaces and whose morphisms are continuous maps.
\end{defi}
\bigskip
The following definition serves to compare pretopologies and is taken from \cite[Section 1.2]{Sta2}:
\begin{defi}
Let $\overline{\bullet}^1$ and $\overline{\bullet}^2$ be two preclosure operators on $X$. We say that $\overline{\bullet}^1$ is finer than $\overline{\bullet}^2$, or $\overline{\bullet}^2$ is coarser than $\overline{\bullet}^1$ if $\overline{A}^1 \subset \overline{A}^2$ for all $A \in P(X)$.
\end{defi}
We can use this definition to define the product pretopology in analogy to how it is usually defined for topological spaces.
\begin{defi}
\label{def:product_preclosure}
Given two pretopological spaces $(X,\overline{\bullet}^1)$ and $(Y,\overline{\bullet}^2)$, the \emph{product preclosure operator} $\overline{\bullet}$ on $X\times Y$ is defined to be the coarsest preclosure operator that makes the canonical projection maps $\pi_1:X\times Y\to X$ and $\pi_2:X\times Y\to Y$ continuous. 
\end{defi}
\begin{lemma}
\label{lem:productpretop}
Given two pretopological spaces $(X,\overline{\bullet}^1)$ and $(Y,\overline{\bullet}^2)$, $\overline{\bullet}$ is the product preclosure operator on $X\times Y$ iff $\overline{B}=\overline{\pi_1(B)}^1\times\overline{\pi_2(B)}^2$ for every subset $B$ of $X\times Y$.
\end{lemma}
\begin{proof}
Continuity of $\pi_1$ and $\pi_2$ requires that $\pi_i(\overline B) \subset \overline{\pi_i(B)}^i, i\in\{1,2\}$. The biggest (and thus coarsest) subset $\overline B\subset X\times Y$ that can fulfill this constraint is the one which fulfills $\pi_i(\overline B) = \overline{\pi_i(B)}^i, i\in\{1,2\}$, which is true iff
\begin{equation}
    \overline B = \{(x,y)~|~x\in \overline{\pi_i(B)}^1,~y\in \overline{\pi_i(B)}^2\} = \overline{\pi_1(B)}^1\times\overline{\pi_2(B)}^2.
\end{equation}
\end{proof}
\begin{coro}
\label{cor:product_in_pretop_cat}
    In the category of pretopological spaces, the product of $(X,\overline{\bullet}^1)$ and $(Y,\overline{\bullet}^2)$ is indeed $X\times Y$ with the product preclosure operator.
\end{coro}
\begin{proof}
The categorical product of $X$ and $Y$ is defined as an object $P$, equipped with two morphisms $\pi_1:P\to X$, $\pi_2:P\to Y$ satisfying the universal property that for every object $Z$, and pair of morphisms $f_1:Z\to X$ and $f_2:Z\to Y$, there exists a unique morphism $f:Z\to P$ such that $f_1=\pi_1\circ f$ and $f_2=\pi_2\circ f$. With $P=X\times Y$ one can always find $f=(f_1,f_2)$ such that in the underlying category of sets $f_i=\pi_i\circ f$.
For $f$ to be continuous for all possible $Z$,
the pretopology on $P$ must be as coarse as possible. 
At the same time, continuity of $\pi_1$ and $\pi_2$ ensure that the coarsest one is the one specified in Lemma \ref{lem:productpretop}.
\end{proof}

The following Lemma now helps us to relate pretopological spaces to cost functions.
\begin{lemma}
  A relation $c$ with \rf{1}, \rf{2} defines a pretopology in $S$ by putting,
  for some $r>0$,
  \begin{equation}\label{25}
  \overline{A}:=\overline{A}_r:=\{ q\in S: c(p,q)\le r \text{ for some }p\in A\}.
  \end{equation}
  A pretopology is also obtained by putting
   \begin{equation}\label{26}
  \overline{A}:=\bigcap_{r>0}\overline{A}_r
  \end{equation}
  with $\overline{A}_r$ as in \rf{25}.
  \end{lemma}
  \begin{proof}
  Hopefully clear.  
  \end{proof}

  Since a cost function $c$ with \eqref{1} and \eqref{2} satisfies all the properties of a metric except symmetry (and possibly finiteness, but that is not relevant here), we can still define closed (open) balls. But now, two types of balls centered at a point  arise,  outward and inward balls:
\begin{align}
B^+(p,t)&:=\{q\in S:\: c(p,q)\leq t\},&\label{outball}\\
B^-(p,t)& :=\{q\in S:\: c(q,p)\leq t\}.& \label{inball}
\end{align}
With the above definitions of outward and inward balls, one can easily see that the preclosure operator (\ref{25}) returns the outward $r$-thickening of a subset $A$, while (\ref{26}) outputs the intersection of such outward thickenings.
We observe
\begin{lemma}
\label{lem:cost_morphisms_are_continuous}
  A cost morphism $f:(S_1,c_1)\to (S_2,c_2)$ is continuous w.r.t.~the pretopologies induced by the cost functions.
\end{lemma}
\begin{proof} Let $r\ge 0$. 
    If $p\in B$ and $q\in \overline{B}$, i.e., $c_1(p,q)\le r$, then $c_2(f(p),f(q))\le r$ since  $f$ is a cost morphism. Hence $f(q)\in \overline{f(B)}$. And by taking limits, this also holds for the pretopology defined in \eqref{26}. 
\end{proof}
\begin{lemma}
    If the pretopologies on $X$ and $Y$ are induced by cost functions $c_1$ and $c_2$, then the product preclosure operator defined in \ref{def:product_preclosure} is induced by $c((x_1,y_1),(x_2,y_2)):=\max(c_1(x_1,x_2),c_2(y_1,y_2))$. With this definition, $(X\times Y,c)$ is the product in the category of cost spaces (cf.~Definition \ref{def:cost_morphisms_and_cost_space_category}).
\end{lemma}
\begin{proof}
    By Lemma \ref{lem:productpretop}, the closed sets in the product pretopology are of the form $\overline B = \{(x,y)~|~x\in \overline{\pi_i(B)}^1\text{ and }y\in \overline{\pi_i(B)}^2\}$. Now, assuming \eqref{25}, $x\in\overline{\pi_i(B)}^1$ iff $\exists b\in \pi_1(B)~:~c_1(b,x)\le r$ and similarly $y\in\overline{\pi_i(B)}^2$ iff $\exists b'\in \pi_2(B)~:~c_2(b',y)\le r$. Hence, both are true simultaneously iff $\exists (b,b')\in B$ such that $\max(c(b,x),c(b',y))\le r$. This proves that the product preclosure is induced by $c$.
    That $(X\times Y,c)$ is the product in the category of cost spaces then follows from Corollary \ref{cor:product_in_pretop_cat} and Lemma \ref{lem:cost_morphisms_are_continuous}.
\end{proof}
\textbf{Remark.}
    The cost function $c$ in the above lemma is in fact the $l_{\infty}$ product of two cost functions $c_1,c_2$. By the lemma, it induces the coarsest topology on $X\times Y$ for which the projections onto the components are continuous. However, it is important to emphasize that one need not restrict oneself to this particular choice. There is a whole spectrum of ways to aggregate the pair $(c_1, c_2)$ , and this flexibility is often crucial in applications. For instance, any $l_p$ product $c_1\times_{l_p}c_2$ for $1\leq p\leq\infty$ defined by $c((x_1,y_1),(x_2,y_2)):=(c_1(x_1,x_2)^p+c_2(y_1,y_2)^p)^{1/p}$ may be chosen as the cost function on $X\times Y$. Yet, other options are  the aggregation via m-schemes \cite{Barth25}, and the coupling of two cost functions \cite{Sturm23}.

\begin{defi}
  A \emph{path} in the pretopological space $X$ is a continuous map
  \bel{top12}
  w:(I,c)\to X
  \qe
  where $c$ is the cost function \rf{2a} and $I$ (here and in the sequel) is equipped with the
  pretopology defined by $c$ (either using (\ref{25}) or (\ref{26}) depending on the application).\\
    A \emph{loop} with base point $p$ is such a continuous map with $w(1)=w(0)=p$. 
\end{defi}
We can use this to study homotopy in the generalized setting of cost functions. To this end we define an equivalence relation on the space of paths using a directed notion of homotopy.
\begin{defi}
\label{def:directed_homotopy}
  The paths $w_1, w_2:I\to X$ with $w_1(0)=w_2(0)=:p, w_1(1)=w_2(1)=:q$ are \emph{equivalent} if there exist continuous maps
  $W_1,W_2:I\times I \to X$ with
  \ba
  W_1(t,0)=w_1(t), & W_1(t,1)=w_2(t),&\\
  W_2(t,0)=w_2(t), &W_2(t,1)=w_1(t)& \text{
    for all }t\in I\\
  W_1(0,s)=W_2(0,s)=p, & W_1(1,s)=W_2(1,s)=q & \text{
    for all }s\in I\ .
  \ea
\end{defi}

\begin{lemma}
\label{lem:surjective_map_is_monotonically_increasing}
    A continuous map $I\to I$ (where, to repeat,  $I$ is equipped with the pretopology defined by the cost function \rf{2a}) is monotonically increasing.
\end{lemma}
\begin{proof}
 For $s\leq t$, considering the pretopology defined by (\ref{25}), 
 we note that, by continuity, $\sigma([s,s+r])=\sigma(\overline{\{s\}})\subset\overline{\sigma(\{s\})}=[\sigma(s),\sigma(s)+r]$. Since $\sigma(s)$ is necessarily mapped to the first point of the interval $[\sigma(s),\sigma(s)+r]$, we can conclude that 
 $\sigma(s)\leq\sigma(t)$ when $s\leq t\leq s+r$. If $t\geq s+r$, then we choose $u_1:= s+\frac{r}{2}$ that implies $\sigma(s)\leq\sigma(u_1)\leq \sigma(s)+r$. If $u_1\leq t\leq u_1+r$, then $\sigma(s)\leq\sigma(u_1)\leq \sigma(t)$. Otherwise we chose $u_2:=u_1+\frac{r}{2}$ and continue until for some $u_i$ we get  $u_i\leq t\leq u_i+r$. 
 \end{proof}
 Lemma \ref{lem:surjective_map_is_monotonically_increasing} in particular  implies that a surjective continuous map $\sigma:I\to I$ preserves the end points, i.e. $\sigma(0)=0$ and $\sigma(1)=1$.
As a consequence, a particular example of equivalent paths is the pair $\left(w ,w\circ \sigma\right)$ when $\sigma:I\to I$ is a surjective continuous map.
We note that if $w$ as in \rf{top12} is continuous, the map $w^\ast$ defined
by $w^\ast(t)=w(1-t)$ need not be continuous, which is why need both $W_1$ and $W_2$ in Definition \ref{def:directed_homotopy}.
\begin{lemma}\label{monoid}
  The equivalence classes of loops with base point $p$ form a monoid $m(p)$.
\end{lemma}
\begin{proof}
  Such loops $w_1,w_2$ can be composed as
  \be\label{eq:pathcomposition}
  w(t):=w_2\circ w_1(t):=\begin{cases} w_1(2t) \text{ for } t\le 1/2\\
    w_2(2t-1) \text{ for } t\ge 1/2.
    \end{cases}
    \qe
    And the constant loop $w_0(t)=p$ for all $t$ is the unit element, since
    for any loop $w$ with base $p$, $w\circ w_0$ and $w_0\circ w$ are
    equivalent to $w$. 
  \end{proof}
This monoid generalizes the fundamental group in homotopy theory to the non-commutative setting adopted here.
  
  We point out that in general, the monoids for different points are not
  necessarily isomorphic. For instance, add to $I$ a loop at the end point
  $1$. Then for $t\in I, t<1$, the monoids $m(t)$ are trivial, because there
  are no loops based at $t$ (continuity prevents them for the pretopology induced by \eqref{2a} on $I$), whereas the monoid 
  at $t=1$, and also each monoid   for every point on the loop that we have added, is isomorphic to $\N$.

  \section{Path spaces and associated Lawvere metric spaces}\label{enrich}

  As already mentioned in Section \ref{weakconditions}, Lawvere \cite{lawvere1973metric} observed the analogy between the triangle inequality 
 \begin{equation}\label{hc1}
     c(q,r) + c(p,q) \ge c(p,r)
 \end{equation}
 and the composition rule in enriched categories
 \begin{equation}\label{hc1a}
    c(q,r) \otimes c(p,q) \to c(p,r).
 \end{equation}
We briefly recapitulate the definition of an enriched category for the reader. 
\begin{defi}\label{enrichdef}
Let  $(V, \otimes, U)$ be a monoidal category with tensor unit $U$ and associator $\alpha:(a\otimes b) \otimes c\to a\otimes (b\otimes c)$.

A $V-$enriched category $S$ consists of: 1) A set of objects $\text{Ob}(S)$; 2) for each ordered pair $(p,q)$ of objects in $S$, an object $c(p,q)$ in  $V$; 3) for each triple $(p,q,r)$ of objects in $S$, a morphism $c(q,r)\otimes c(p,q)\to c(p,r)$ in $V$; and 4) for each object $p$ in $S$, a morphism $U \to c(p,p)$ in $V$; such that $\forall a,b,c,d \in \text{Ob}(S)$:
\begin{itemize}
    \item Composition is associative:
    \begin{center}
        \begin{tikzcd}
        (c(c,d)\otimes c(b,c))\otimes c(a,b) \arrow[rr, "\alpha"] \arrow[d] &          & c(c,d)\otimes (c(b,c)\otimes c(a,b)) \arrow[d] \\
        c(b,d)\otimes c(a,b) \arrow[r]                                      & c(a,d) & c(c,d)\otimes c(a,c) \arrow[l]
        \end{tikzcd}
    \end{center}
    \item Composition is unital:
    \begin{center}
        \begin{tikzcd}
{c(b,b)\otimes c(a,b)} \arrow[r]       & {c(a,b)} & {c(a,b)\otimes c(b,b)} \arrow[l]       \\
{U\otimes c(a,b)} \arrow[u] \arrow[ru] &          & {c(a,b)\otimes U} \arrow[u] \arrow[lu]
\end{tikzcd}
    \end{center}
\end{itemize}

\end{defi}

For the special case in which $(V,\otimes, U)$ is a poset of numbers with addition as tensor product $((P,\ge),+,0)$, the composition $c(q,r)\otimes c(p,q)\to c(p,r)$ then indeed turns into the triangle inequality $c(q,r) + c(p,q)\ge c(p,r)$.
That $c(b,b) = 0$ for all $b$ follows from the fact that composition is unital: We have $c(a,b)=0+c(a,b)\ge c(b,b)+c(a,b)\ge c(a,b)$ and thus for $a=b$ we get $c(b,b)\ge 2c(b,b)\ge c(b,b)$. For $P=[0,\infty]$, we obtain a Lawvere metric space. 

\bigskip

Given a pretopological space $X$ with a cost function $c$, we next define natural categories for its paths and look at associated Lawvere metric spaces. In fact, there are two different types of relations between paths that we can use for  defining categories.
In the first case, we define a category which is like a directed version of a fundamental groupoid: Objects are points of $X$ and morphisms are directed homotopy classes of paths, where
directed homotopies are defined as in \ref{def:directed_homotopy}. This gives rise to a groupoid when every path $\gamma$ has an inverse $\gamma^{-1}(t):=\gamma(1-t)$ such that the composition $\gamma\circ \gamma^{-1} \sim \text{id}$ 
where $\text{id}$ is the constant path at a point. One could then consider homotopies of homotopies to obtain a directed infinity groupoid, generalizing the infinity groupoid obtained from homotopies between paths of topological spaces.\\
\medskip

The second way for constructing such a category  utilizes the cost function. 
We first introduce the standard definition of the length of a continuous curve (or path, we shall use the two terms synonymously) $\gamma:[a,b]\to X$ \cite{Sche}, 
\begin{equation}\label{le1}
   \lth(\gamma):=\sup_{a=t_0<t_1<t_2 <\dots <t_m=b} \sum_{\mu =1}^m c(\gamma(t_{\mu -1}),\gamma(t_\mu))
\end{equation}
where the supremum is taken over all partitions of $[a,b]$. A curve with finite length is called rectifiable, and in the sequel we consider the class of rectifiable curves. The length of a curve is invariant under reparametrization. That allows us to select a particular parametrization. A rectifiable curve can be parametrized by arclength, i.e., on the interval $[0,b]$ with fixed $a=0$ and
\begin{equation}\label{le2}
   t= \lth(\gamma_{[0,t]})
\end{equation}
and $b=\lth(\gamma)$. This is obtained from an arbitrary parametrization by inverting the function $\ell(\tau)=\lth(\gamma_{[0,\lth(\gamma_{[0,\tau]})]})$. In particular, a constant curve then is parametrized on the point interval $[0,0]$. \\
If we want to have all curves parametrized on the same interval, which we can take as $[0,\infty)$, we may put
\begin{equation*}
    \gamma(t)=\gamma(b) \text{ for }t\ge b
\end{equation*}
so that it becomes constant on $[b,\infty)$. \\
Alternatively, we can parametrize $\gamma:[0,b]\to X$ proportionally to arclength on $[0,1]$ by
\begin{equation}\label{le3}
   s= \frac{\lth(\gamma_{[0,s]})}{\lth(\gamma_{[0,b]})}.
\end{equation}
That has the advantage that all (rectifiable) curves are parametrized on the same interval $[0,1]$.\\

\medskip

\begin{rem}
    In contrast to the setting of Moore path spaces \cite{BR} where a length function is simply assumed on a topological space, we here work with a length function that comes from a cost function (a metric in the symmetric case). We can then parametrize paths proportionally to or by arclength so that we can naturally restrict to subpaths.
\end{rem}

We shall now introduce two categories that we can enrich in the sequel.
The first category is $S_X$   whose objects are points of $X$ and whose morphisms are  rectifiable paths between two given points, parametrized either by arclength as in \eqref{le1} or proportionally to arclength as in \eqref{le3}. For the composition of paths  with \eqref{le1}, let $\gamma_1:[0,b_1]\to X, \gamma_2:[0,b_2]\to X$ with $\gamma_{2}(0)=\gamma_{1}(b_1)$ and put
\begin{eqnarray}
\nonumber
\gamma_{12}=\gamma_2\circ \gamma_1&:&[0,b_1+b_2]\to X\\
\gamma_{12}(t)&=&\begin{cases}
    \gamma_{1}(t)& \text{ for } 0\le b_1\\
    \gamma_{2}(t-b_1) & \text{ for } b_1\le t\le b_2
\end{cases}
\label{le6}
    \end{eqnarray}
    and extend it as a constant curve for $t\ge b_1+b_2$.\\
    With \eqref{le3} and $\gamma_1, \gamma_2:[0,1]\to X$ with $\gamma_{2}(0)=\gamma_{1}(1)$, we put
\begin{eqnarray}
\nonumber
\gamma_{12}=\gamma_2\circ \gamma_1&:&[0,1]\to X\\
\nonumber
\gamma_{12}(t)&=&\begin{cases}
    \gamma_{1}(\frac{t}{\ell})& \text{ for } 0\le t\le \ell\\
    \gamma_{2}(\frac{t-\ell}{1-\ell}) & \text{ for } \ell \le t\le 1
\end{cases}\\
&&\text{with } \ell =\frac{\lth(\gamma_1)}{\lth(\gamma_1) +\lth(\gamma_2)}
\label{le7}
    \end{eqnarray}
Since for a constant curve $\gamma_0$, $\lth(\gamma_0)=0$, in either case, existence of identities and associativity of composition is guaranteed. Thus, both constructions define a category.\\
The second category is $T_X$ whose objects again are points $p,q$ of $X$ and where we put a morphism $(p,q)$ whenever there is a rectifiable path from $p$ to $q$. Since the constant path is rectifiable and the composition of rectifiable paths is again rectifiable, $T_X$ is a category, indeed. 
\medskip

From these categories we will now construct  monoidally enriched categories.\\
In order to construct an enriched category related to category $S_X$, 
we introduce the monoidal category $V$  with objects being sets of non-negative numbers, morphisms $\subseteq$, and $\otimes=+$ as the summation of numbers in a pair of sets
\[
A\otimes B=\{a+b|a\in A,b\in B\}
\]
with the identity $\{0\}$.
The  category $S_X$ with $X$ the set of objects and rectifiable paths as hom-set $\Hom_{S_X}(p,q)$ then gives rise to an  enriched category $P_X$, in which the objects are those of $S_X$, while the morphisms are lengths of paths. Then we can define a functor $\pi_S$ from $S_X$ to that enriched category that is the identity on objects and applies the $\lth$-operator to paths, 
\begin{equation}\label{pi1}
    \pi_S(\gamma)=\text{length}(\gamma).
    \end{equation}
 $\pi_S$ sends identities (constant paths) to identities (the tensor unit $0$) and, whenever we have  paths $\gamma_1$ from $p$ to $q$ and $\gamma_2$ from $q$ to $r$, then we have 
     \begin{equation}\label{hc20}
     \lth(\gamma_2\circ \gamma_1)=\lth(\gamma_1) +\lth(\gamma_2).
     \end{equation}
   Thus, $\pi_S$ is indeed a functor.

Another  Lawvere metric space $Q_X$ emerges when we take as objects  again the points of $X$ and define the  morphisms by 
     \begin{equation}\label{hc21}
     c(p,q): = \inf_{\gamma \in S_X(p,q)} \lth(\gamma)
     \end{equation}
     where $S_X(p,q)$ is the space of paths from $p$ to $q$. 
     With this definition, \eqref{hc1} is again ensured, making $Q_X$ into a Lawvere metric space. There is a functor $\xi:P_X \to Q_X$, that acts as the identity on objects. For the action on morphisms, we define 
     \begin{equation}
        \xi(L)=\inf_{L'}\{~L'~|~\text{dom}(L)=\text{dom}(L'),~\text{codom}(L)=\text{codom}(L')~\}
     \end{equation}
     or 
     \begin{equation}
        \xi(L \in P_X(p,q))=\inf_{L'\in P_X(p,q)}L'.
     \end{equation}
     But for the present  construction,  we can also work with more general length concepts than above. \\
     
     A curve $\gamma\in S_X(p,q)$ with
     \begin{equation}
         \lth(\gamma)=\pi_S(\gamma)=\xi\circ \pi_S(\gamma)=c(p,q)
     \end{equation}
     then is a tachistic in the sense of Def. \ref{tach}. A chronodesic in the sense of Def. \ref{general} then satisfies the length minimizing property for any two  points on them that are connected by sufficiently short subpaths. \\  
      We can also define the betweenness relation with this construction:
For every $p,q,r \in Q_X$,  $b(p,q,r)$ holds if and only if composition admits a reverse arrow, that is, there is an arrow $c(p,r)\to c(q,r) \otimes c(p,q)$ (making the triangle inequality into an equality). 
On a tachistic $\gamma:[a,b]\to X$ from $p$ to $r$, every other point $q$ then is between $p$ and $r$.\\
\medskip

   \section{Curvature}\label{curv}  

The notion of total convexity can be generalized in the following way: \\
$S$ is totally convex, if for every pair $p,r\in S$ (where $S$ is equipped with a cost function $c$) and for every two objects $t_1,t_2\in P  =(([0,\infty],\ge),+)$ such that there is a morphism $t_1+ t_2\to c(p,r)$,  there exists $q\in S$ such that there are morphisms $t_1\to  c(p,q)$ and $t_2\to c(q,r)$. In the context of a length structure, this implies that any tachistic
with length $\leq t_1+t_2$ can be decomposed (by the means of concatenation) into two tachistics $\gamma_1$ and $\gamma_2$ with lengths $\leq t_1$ and $\leq t_2$ respectively.

The property of a space being almost chronodesic (or of a path being tachistic) can also be defined: $S$ is almost chronodesic if, for every  $p,r\in S$ and $\epsilon >0$, whenever $t_1+t_2=c(p,r)$, then there exists $q\in S$ s.t. $t_1+\epsilon\to c(p,q)$ and $t_2+\epsilon\to c(q,r)$.\\

The concepts of median and triple betweenness naturally extend the betweenness relation $b(p,q,r)$. For a triple  $(x_1,x_2,x_3)$, finding an intermediate point requires selecting a direction for each pair and defining the median as a point between every pair relative to these directions. For instance, one can say $b(x_1,x_2,x_3:x_\star)$ holds if and only if $b(x_1,x_\star,x_2)$, $b(x_2,x_\star,x_3)$ and $b(x_3,x_\star,x_1)$  all hold.  And $x_\star$ is then called a \emph{median} of the directed triple $(x_1,x_2,x_3)$.

In  terms of the balls \eqref{outball} and \eqref{inball}, the conditions for a median can  be expressed as follows:
\begin{align}\label{3inter}
x_\star&\in B^+(x_1,r_1)\cap B^-(x_2,r'_2)&\\
x_\star&\in B^+(x_2,r_2)\cap B^-(x_3,r'_3)&\\
x_\star&\in B^+(x_3,r_3)\cap B^-(x_1,r'_1)&
\end{align}
 where 
\[
r_i:=c(x_i,x_\star),\quad r'_i:=c(x_\star,x_i) \text{ for } i=1,2,3.
\]

Admitting a median is a property that depends on the prior ordering of the points and is not invariant under permutation.\\

In general, while the inequalities $r_1+r'_3\ge c(x_1,x_3)$ etc. are satisfied, there may still not exist a point $x_\star$ that satisfies all the betweenness conditions, or equivalently, satisfies \eqref{3inter}. 

We can then enlarge the radii of the balls until we obtain a common point of intersection. The scaling factor required to achieve this intersection provides a measure of the deviation from admitting a median.

  When we want to achieve some generalization of
our curvature notions for metric spaces, it seems important to take as the 
basic elements not the points, but rather the directed pairs $(p,q)=p\to
q$. We have called $t$ a median of $(p,q), (q,r), (r,p)$ if
\bel{10}
c(p,q)=c(p,t)+c(t,q),\ c(q,r)=c(q,t)+c(t,r),\ c(r,p)=c(r,t)+c(t,p).
\qe
This is intended as a non-symmetric analogue of a metric tripod. There one has a distance function $d(.,.)$ and three points $p,q,r$ and a median $t$ with the property that $d(p,q)=d(p,t)+d(t,q), d(q,r)=d(q,t)+d(t,r), d(r,p)=d(r,t)+d(t,p)$. In other words, in that situation, we have an undirected graph with vertices $p,q,r,t$ and edges connecting each of $p,q,r$ with $t$. \\
Of course, we can do the same as in \eqref{10} for any three directed pairs involving three
points. And we can require analogues for more than three pairs. Requiring
\rf{10} for any collection of directed pairs would constitute an analogue of
hyperconvexity. And curvature would again quantify the minimal deviation from
\rf{10}, that is, to what extent the best possible choice of $t$ would violate
the equalities in \rf{10}. \\

We could thus define directed curvature $\rho$ in analogy to the undirected case, given in \cite{joharinadJost2019}, by the following equation:
\begin{equation}
  \label{eq:directedCurvature}
  \begin{split}
    \rho(x_1,x_2,x_3) = \sup_{r_1,r_2,r_3\ge 0} \left\{ \inf_x \max_{i\in \{1,2,3\}} \frac{c(x_i,x)}{r_i} ~\bigg|~
    \begin{pmatrix}
      r_1+r_2\ge c(x_1,x_2)\\
      r_2+r_3\ge c(x_2,x_3)\\
      r_3+r_1\ge c(x_3,x_1)
    \end{pmatrix}
     \right\}
  \end{split}
\end{equation}
In contrast to the symmetric case, care must be taken in the specification of the constraints $r_i+r_j \ge c(x_i,x_j)$ because $c$ is directed.
As in the undirected case, the sup-inf is achieved exactly when $r_i+r_j=c(x_i,x_j)$, and we obtain a system of 3 equations with 3 unknowns from which the so-called Gromov products can be computed (again, the direction of $c$ has to be respected in the process):
\begin{equation}
  \begin{split}
    2r_1 &= (r_1+r_2)+(r_3+r_1)-(r_2+r_3) = c(x_1,x_2)+c(x_3,x_1)-c(x_2,x_3),\\
    2r_2 &= (r_1+r_2)+(r_2+r_3)-(r_3+r_1) = c(x_1,x_2)+c(x_2,x_3)-c(x_3,x_1),\\
    2r_3 &= (r_2+r_3)+(r_3+r_1)-(r_1+r_2) = c(x_2,x_3)+c(x_3,x_1)-c(x_1,x_2).
  \end{split}
  \label{eq:directedCurvatureGromovIdentities}
\end{equation}
In contrast to the symmetric case $\rho(x_1,x_2,x_3)$ has fewer symmetries and for every of the $3!$ different ways to order $x_1,x_2,x_3$ we obtain possibly different values.

\begin{rem}
Alternatively, one can also define a symmetrized version of curvature. To this end, note that any non-symmetric cost function or Lawvere metric $c$ can be symmetrized with any symmetric binary operator that respects the directed triangle inequality. For example, a canonical choice is 
\begin{equation}
  \begin{split}
    d(x,y):=\frac{1}{2}(c(x,y)+c(y,x)).
     \end{split}
\end{equation}
Having defined $d$, one can then define the symmetrized curvature using \eqref{eq:directedCurvature} but with $d$ in place of $c$. Depending on which aspect of the space under consideration one is interested in, this can already yield interesting information.
\end{rem}

\begin{rem}
    Combining our definition of directed curvature \eqref{eq:directedCurvature} with our definition of $\mathcal{O}$-cost functions \rf{asymptoticCosts}, we obtain the possibility to speak about the computational geometry of formal languages.
    \end{rem}

\section{Hyperconvex hulls}\label{hyp}
In this section, we shall propose a non-symmetric version of the construction of the tight span in \cite{D}. For a set $(S,c)$ with a finite cost function, let us consider all function pairs
$(f,g)$ with 
\begin{equation}
  \label{hyp1}
  f(p)=\sup_{q\in S}\ (c(p,q)-g(q)).
\end{equation}
In the symmetric case, we may put $f=g$. More precisely, 
if $c$ is symmetric, then those functions $f$,
for which the pairs $(f,f)$ satisfy
\eqref{hyp1}, define the tight span of the
metric space $(S,d)$ in the sense of \cite{D}, or the hyperconvex hull of \cite{Isb}.\\
To illustrate the concept in the non-symmetric case, we recall the example of the unit interval $I=[0,1]$ equipped with the
cost function \rf{hyp2}
\begin{equation*}
  c_1(p,q)=
  \begin{cases}
    q-p& \text{ if }p\le q\\
    2(p-q)& \text{ if }p\ge q.
  \end{cases}
\end{equation*}
For $0\le \alpha \le 1$, we define
\begin{eqnarray}
  f(p)&=
  \begin{cases}
    \alpha -p& \text{ if }p\le \alpha\\
    2(p-\alpha)& \text{ if }p\ge \alpha
  \end{cases}\\
  g(q)&=
  \begin{cases}
    2(\alpha -q)& \text{ if }q\le \alpha\\
    q-\alpha& \text{ if }q\ge \alpha.
  \end{cases}
\end{eqnarray}
Then \rf{hyp1} holds. In fact, we could have started with the endpoints $0,1$
of $I$ and restricted the cost function $c_1$ to them. Again, \rf{hyp1} would
have yielded the above pair $(f,g)$.

\end{document}